\documentclass[10pt]{article}
\usepackage[latin1]{inputenc}
\usepackage{amsmath,amsfonts,amssymb}
\usepackage{amscd}
\usepackage{amscd,graphicx}

\newcommand{\qed}{{$\diamond$ }}

\begin{document}

\title{A theorem on the cores of partitions} 

\author{
J\o rn B. Olsson\\
\small Department of Mathematical Sciences, University of Copenhagen \\[-0.8ex]
\small Universitetsparken 5,DK-2100 Copenhagen \O, Denmark
}

\maketitle

\bigskip

Suppose that $s,t \in \mathbb{N}$ are relatively prime positive
integers. In the study of block inclusions between $s$- and $t$-blocks
of partitions \cite{OS} we introduced an $(s,t)$-abacus to study 
relations between $s$- and $t$-cores of partitions. 
This is because the cores determine the blocks. 

Before we state the main result of this paper 
let us mention that the basic facts about
partitions, hooks and blocks of partitions may be found in \cite{JK}, 
Chapter 2 or \cite{Ol}, Chapter 1.  You may get to the $s$-core
$\lambda_{(s)}$ of a partition $\lambda$ by removing a series of
$s$-hooks (ie. hooks of length $s$) until all $s$-hooks are removed. 
The $s$-core is independent of the order in which the $s$-hooks are
removed. A partition has by definition {\it $s$-weight} $w,$ if you need to
remove exactly $w$ $s$-hooks to get to its $s$-core. It also equals the
number of hooks in the partition of length divisible by $s$
(\cite{JK}, 2.7.40). Thus a partition is an $s$-core if and only if it
has $s$-weight 0. Two partitions of $n$ are said to be in the same
$s$-block if they have the same $s$-core. This definition is of
course inspired by the theorem about irreducible characters of the
symmetric groups, which is still referred to as the Nakayama
conjecture. See \cite{JK}, 6.2.21. The hook structure of a partition
is conveniently determined by its {\it first column hook lengths}, or
more generally any of its {\it $\beta$-sets}
\cite{Ol}, section 1.   
 
In this note we want to illustrate the usefulness of the
$(s,t)$-abacus by showing the following result:

\medskip

\noindent {\bf Theorem 1:} {\it Let  $s,t$ be relatively prime positive integers.
Suppose that $\rho$ is a $t$-core. Then the $s$-core of $\rho$ is again a $t$-core. } 
%\end{theorem}

\medskip

For a non-expert this may not seem to be surprising. However examples
show that in  
removing $s$-hooks from the partition $\rho$ of the theorem you may
have to go through arbitrarily long series of intermediate partitions 
which cannot be chosen as $t$-cores, and the overall behaviour may
indeed appear
to be rather chaotic. Yet the final result turns out to
be again a $t$-core. The proof of Theorem 1  is surprisingly simple, once you
understand how to use the $(s,t)$-abacus. The Theorem plays an
important role in a forthcoming paper \cite{O2} 

\medskip

As we shall see below there is an analogous result for bar partitions
and bar cores (Theorem 4) under the additional assumption that 
$s$ and $t$ are both odd. The proof in this case is somewhat more
delicate. 

\medskip
 
A partition which is at the same time an $s$- and $t$-core is called 
an $(s,t)${\it -core}. 
There are only finity many such partitions and the maximum one has
cardinality  $m_{s,t}=(s^2-1)(t^2-1)/24.$ (\cite{Anderson-pq},
\cite{OS}). 
Therefore our result implies the following:

\smallskip

\noindent {\bf Corollary 2:} {\it A $t$-core of $n$ has $s$-weight at least
$w=\lfloor(n-m_{s,t})/s \rfloor.$ 
Thus it contains at least $w$ hooks of length divisible by $s.$
}

\smallskip

As an example, $m_{3,5}=8$ so that the 3-core $(8,6,4,2^2,1^2)$ of 24
must contain at least 3 hooks of length divisible by 5. 
It actually contains 4 such hooks.

\smallskip

There are a few obvious questions you might ask after having seen the
theorem, but unfortunately they seem to have negative answers:

If $\lambda$ is a partition define
$$\lambda_ {(s,t)}=(\lambda_{(s)})_{(t)}.$$

Thus  $\lambda_{(s,t)}$ is the $t$-core of the $s$-core of $\lambda.$ 
By our theorem  $\lambda_ {(s,t)}$ is actually an $(s,t)$-core and we
may call it the  $(s,t)$-core of $\lambda.$ 
But {\it generally}  $\lambda_ {(s,t)} \neq \lambda_{(t,s)}.$ 
Indeed, if for example 
$\lambda=(3), s=2, t=3$ then $\lambda_ {(2,3)}=(1)$ 
whereas  $\lambda_ {(3,2)}$ is the empty partition. 
 
Also we may define a $(s,t)$-block of $n$ as the set of partitions with the
same $(s,t)$-core. Obviously it is a union of $s$-blocks of $n.$ But
it is not necessarily a union of $t$-blocks. Indeed the partitions of 5 with
empty $(2,3)$-core $(0)$ are $(4,1)$ and $(2,1^3).$ They form a
2-block (of weight 1) but obviously not a 3-block.

Finally, it is known that the number of partitions of $n$ with a 
given $s$-core only depends on the $s$-weight, (\cite{JK}, 2.7.17).
But it is not true, that the number of $t$-cores of $n$ with a given $s$-core
only depend on the $s$-weight.  Indeed here are some examples:

Weight 1:  
The number of 5-cores of 10 with 7-core (3) is 2 
and the number of 5-cores of 10 with 7-core (2,1) is 4. 

Weight 3:
The number of 5-cores of 10 with 3-core (1) is 8.
The number of 5-cores of 11 with 3-core (2) is 3.

\bigskip

Theorem 1 is proved below. As an application we list the following
result. Define the {\it principal} $s$-block of $n$ to be the $s$-block
containing the partition $(n).$ 

\medskip

\noindent {\bf Corollary 3:} {\it Assume that $s>r\geq t$ and that $s$
  and $t$ are relatively prime. Let $n=as+r,$ $a \ge 0$. 
Then any $t$-core of $n$  is {\rm not} contained in the principal $s$-block of $n$. }

\medskip

\noindent{\bf Proof:} If $\rho$ is a $t$-core of $n,$ then its
$s$-core $\rho'$ is a partition of $m=bs+r$ for some $b,$ $a \ge b \ge
0.$ By Theorem 1 $\rho'$ is also a $t$-core. If $\rho$ were in the
principal $s$-block of $n$ then $\rho'=(m).$ But $m \ge r \ge t,$  so
that $(m)$ contains a $t$-hook and thus is not a $t$-core.  
\qed

\medskip

There is an analogous result to Theorem 1 for bar partitions 
under the additional
assumption that $s$ and $t$ are odd. A {\it bar partition} is a
partition into distinct parts. For these partitions there is for odd
integers a theory of {\it bars} corresponding to the theory of hooks
in arbitrary partitons. See eg. \cite{Ol}, Section 4 for details. In particular each
bar partition has for any given odd integer $s$ a unique  $\bar s$-core 
obtaind by removing a series of $s$-bars from the partition. We have

\medskip

\noindent {\bf Theorem 4:} {\it Let  $s,t$ be relatively prime odd positive integers.
Suppose that $\rho$ is a $\bar t$-core. Then the $\bar s$-core of
$\rho$ is again a $\bar t$-core.}

\medskip

There is an analogue of Corollary 2 for bar partitions, 
but the statement is less precise.

\medskip

We assume now that $s,t \in \mathbb{N}$ are relatively prime positive
integers. 

The $s$-abacus was introduced by G. James. Its relation to the study
of $s$-cores and $s$-quotients of partitions is explained in detail in 
\cite{JK}, Section 2.7. (Or see Section 3 in \cite{Ol}.) The  $s$-abacus  has
$s$ infinite runners, numbered $0,...,(s-1)$ going from north to
south. The $i$'th runner contains the nonnegative integers which are
congruent to $i$ modulo $s$ in increasing order. Here is part of the
$7$-abacus:

{\scriptsize 
\begin{center}

\begin{tabular}{l c c c c c c c}
 Runner: &0&1&2&3&4&5&6 \\
\hline
\hline
 &0&1&2&3&4&5&6 \\
 &7&8&9&10&11&12&13 \\ 
 &14&15&16&17&18&19&20 \\ 
 &21&22&23&24&25&26&27 \\
 &28&29&30&31&32&33&34 \\
  &.... \\ 
\end{tabular}
\end{center}
}

Generally we may arrange the first column hook lengths of the maximum 
$(s,t)$-core $\kappa_{s,t}$ in a diagram, which we call the
$(s,t)$-{\it diagram} \cite{Anderson-pq},\cite{OS}.
  
Start with the largest entry $st-s-t$
in the lower left hand corner and subtract multiples of $s$ along the rows and
multiples of $t$ along the columns as long as possible. 
Then the first 
column hooklengths of any $(s,t)$-core must be among the numbers of this
diagram. The reason is, that  $st-s-t$ is the largest integer which
cannot be written in the form $as+bt$ where $a,b$ are non-negative
integers. More details may for example be found in   \cite{Anderson-pq}.

Here is the (5,7)-{\it diagram.}

{\scriptsize
\begin{center}

\begin{tabular}{c c c c c}
2\\ 
9&4\\
16&11&6&1\\
23&18&13&8&3\\
%\hline
\end{tabular}
\end{center}
}

Note that the numbers in the 
columns (read from north to south) are part of some runners on the
$t$-abacus. In the example we have runners 1,2,3,4,6 of the 7-abacus represented. 
The order of the runners is changed and some 
runners are missing. The missing runners may be represented by
extending the diagram to the
south like this:

{\scriptsize
\begin{center}

\begin{tabular}{c c c c c c c}
2\\ 
9&4\\
16&11&6&1\\
23&18&13&8&3\\
%\hline
{\it 30} & {\it 25}& {\it 20}&{\it 15}&{\it 10}&{\it 5}&{\it 0}
\end{tabular}
\end{center}
}

We have just added part of the 0'th  5-runner to the south 
written in {\it italics}. In this diagram also the missing 7-runners  
numbered 0 and 5 are represented.

If we continue adding rows to the south we get the $t$-abacus with
runners in a different order, with numbers increasing from north to south. 
We refer to this as the {\it  $(s,t)$-abacus}.

 The  $(s,t)$-abacus and a numbering of its rows is illustrated by
 this example ($s=5, ~t=7)$:

\bigskip
{\scriptsize
\begin{center}
\begin{tabular}{l c c c c c c c c}

Runner: & 2&4&6&1&3&5&0 \\
\hline
\hline
Row -4&2\\
Row -3&9&4\\ 
Row -2&16&11&6&1\\
Row -1&23&18&13&8&3\\
\hline
%\end{tabular}\\
%\begin{tabular}{c c c c c }
Row 0&30&25&20&15&10&5&0\\
Row 1&37&32&27&22&17&12&7\\
Row 2&44&39&34&29&24&19&14\\
Row 3&51&46&41&36&31&26&21\\
Row 4&58&53&48&43&38&33&28\\
\hline
Row 5&65&60&55&50&45&40&35\\
Row 6&72&67&62&57&52&47&42\\
Row 7&79&74&69&64&59&54&49\\
Row 8&86&81&76&71&66&61&56\\
....
\end{tabular}
\end{center}
}
\bigskip

The rows in the $(s,t)$-diagram are numbered -1,-2,..., starting from
the bottom. 
The rows below the $(s,t)$-diagram are numbered 0,1,2... starting from
the top 
as indicated in the example. Thus the $i$-th row contains a decreasing
sequence of numbers which are congruent modulo $s$. The difference
between neighbouring numbers is $s$ and the eastmost number in the
row (on the runner $0$) is $t\cdot i.$ Clearly any non-negative integer is
represented uniquely by a position on the $(s,t)$-abacus. The runners
of the $s$-abacus are visible in the rows of the
$(s,t)$-abacus. Rows, whose numbers differ by a multiple of $s$ (like
rows -2,3 and 8 in the example) contain numbers from the same runner
of the $s$-abacus. Thus the $s$-runners are broken into pieces. 

\medskip

This means that (as already mentioned) the $(s,t)$-abacus is useful 
for studies involving the relations
between $s$-cores  and $t$-cores. Take a $\beta$-set for a given
partition and represent its numbers as beads on the  $(s,t)$-abacus. 
This means that we place a bead in the position numbered $i$ on the
abacus for all $i$ in the $\beta$-set. Adding/removing $s$-hooks from partitions are
reflected by horisontal moves of the beads. Adding/removing $t$-hooks are
reflected by vertical moves of the beads on the $(s,t)$-abacus where
a ``vertical move'' could include a shift of $s$ rows, corresponding
to the breakup of the $s$-runners we have just described. 

We are now in the position to prove Theorem 1. 
Of course a $t$-core need not be an $s$-core. We show that the $s$-core of an
$t$-core is an $(s,t)$-core.

\medskip

\noindent{\bf Proof of Theorem 1:} Suppose that $\rho$ is a $t$-core. 
Let $X$ be the set of first column hook lengths of $\rho.$ 
We extend $X$ to a larger $\beta$-set $Y$ for $\rho$ in such a way that it contains 
all the numbers in the $(s,t)$-diagram. We represent the numbers in $Y$ as 
beads on the $(s,t)$-abacus. 
The result is a diagram where there is no empty space 
to the north of a bead (since $\rho$ is a $t$-core) and the part consisting 
of the $(s,t)$-diagram is filled with beads. Each runner contains a number of beads 
outside the $(s,t)$-diagram. 
If row $i \ge 0$ contains 
$n_i \ge 0$ beads, then clearly $n_0 \ge n_1 \ge n_2 \ge ... \ge n_i \ge ...$ 

The removal of an $s$-hook from $\rho$ is reflected by moving a bead
to an 
empty space next to it to the east or (if it on the eastmost runner) to an 
empty space at the westmost runner $s$ rows above. 
We have reached a bead configuration for the $s$-core of $\rho$ when no 
more moves of this kind are possible. 

We adapt the following strategy for the moves.
 
{\it Step 1:} Start by moving all beads as far to the east as possible 
in their respective rows. Then the number of beads outside the 
$(s,t)$-diagram on the runners is decreasing, when we move from the
west to the 
east, by the remark on the $n_i$'s above.  Moreover the beads 
still represent a $t$-core.    

{\it Step 2:} Move all beads on runner 0 to the westmost runner in the section above it, 
($s$ rows above it) if there is an empty space. The beads still represent a $t$-core.

Then repeat step 1 and 2 for as long as possible. The process stops when 
we have reached a 
bead configuration for the $s$-core of  $\rho.$ Since each step results in a $t$-core, 
the result follows. \qed

\medskip

Here is an example illustrating the steps of the proof ($s=5, t=7$): 
The numbers in boldface in the first diagram are the ones in the 
$\beta$-set $Y,$ as described in the
proof. We have that the initial 7-core is
$$\rho=(42,36,30,24,18,12,11,7,6,2^4,1),$$
 a partition of 195. 
We apply first step 1, then step 2 and then step 1
again. Rows without beads are omitted.
The 5-core of $\rho$ is $(5,4,2,1^4)$, another 7-core (of 15), and the
5-weight of $\rho$ is 36.

{\scriptsize
\begin{center}
\begin{tabular}{l c c c c c c c c}

Runner: & 2&4&6&1&3&5&0 \\
\hline
\hline
Row -4&{\bf 2}\\
Row -3&{\bf 9}&{\bf 4}\\ 
Row -2&{\bf 16}&{\bf 11}&{\bf 6}&{\bf 1}\\
Row -1&{\bf 23}&{\bf 18}&{\bf 13}&{\bf 8}&{\bf 3}\\
\hline
%\end{tabular}\\
%\begin{tabular}{c c c c c }
Row 0&{\bf 30}&{\bf 25}&20&{\bf 15}&{\bf 10}&{\bf 5}&{\bf 0}\\
Row 1&37&{\bf 32}&27&22&{\bf 17}&12&{\bf 7}\\
Row 2&44&{\bf 39}&34&29&24&19&14\\
Row 3&51&{\bf 46}&41&36&31&26&21\\
Row 4&58&{\bf 53}&48&43&38&33&28\\
\hline
Row 5&65&{\bf 60}&55&50&45&40&35\\
Row 6&72&{\bf 67}&62&57&52&47&42\\
%Row 7&79&{\bf 74}&69&64&59&54&49\\
%Row 8&86&{\bf 81}&76&71&66&61&56\\
....
\end{tabular}

\end{center}
}

{\scriptsize
\begin{center}
\begin{tabular}{l c c c c c c c c}

Runner: & 2&4&6&1&3&5&0 \\
\hline
\hline
Row -4&{\bf 2}\\
Row -3&{\bf 9}&{\bf 4}\\ 
Row -2&{\bf 16}&{\bf 11}&{\bf 6}&{\bf 1}\\
Row -1&{\bf 23}&{\bf 18}&{\bf 13}&{\bf 8}&{\bf 3}\\
\hline
%\end{tabular}\\
%\begin{tabular}{c c c c c }
Row 0&30&{\bf 25}&{\bf 20}&{\bf 15}&{\bf 10}&{\bf 5}&{\bf 0}\\
Row 1&37&32&27&22&{\bf 17}&{\bf 12}&{\bf 7}\\
Row 2&44&39&34&29&24&19&{\bf 14}\\
Row 3&51&46&41&36&31&26&{\bf 21}\\
Row 4&58&53&48&43&38&33&{\bf 28}\\
\hline
Row 5&65&60&55&50&45&40&{\bf 35}\\
Row 6&72&67&62&57&52&47&{\bf 42}\\
%Row 7&79&74&69&64&59&54&{\bf 49}\\
%Row 8&86&81&76&71&66&61&{\bf 56}\\
....
\end{tabular}
\end{center}
}

{\scriptsize
\begin{center}
\begin{tabular}{l c c c c c c c c}

Runner: & 2&4&6&1&3&5&0 \\
\hline
\hline
Row -4&{\bf 2}\\
Row -3&{\bf 9}&{\bf 4}\\ 
Row -2&{\bf 16}&{\bf 11}&{\bf 6}&{\bf 1}\\
Row -1&{\bf 23}&{\bf 18}&{\bf 13}&{\bf 8}&{\bf 3}\\
\hline
%\end{tabular}\\
%\begin{tabular}{c c c c c }\hline
%\end{tabular}\\
%\begin{tabular}{c c c c c }
Row 0&{\bf 30}&{\bf 25}&{\bf 20}&{\bf 15}&{\bf 10}&{\bf 5}&{\bf 0}\\
Row 1&{\bf 37}&32&27&22&{\bf 17}&{\bf 12}&{\bf 7}\\
Row 2&44&39&34&29&24&19&{\bf 14}\\
Row 3& 51&46&41&36&31&26&{\bf 21}\\
Row 4&58&53&48&43&38&33&{\bf 28}\\
%\hline
%Row 5&65&60&55&50&45&40&35\\
%Row 6&72&67&62&57&52&47& 42\\
%Row 7&79&74&69&64&59&54& 49\\
%Row 8&86&81&76&71&66&61&56\\
....
\end{tabular}
\end{center}
}

{\scriptsize 
\begin{center}
\begin{tabular}{l c c c c c c c c}

Runner: & 2&4&6&1&3&5&0 \\
\hline
\hline
Row -4&{\bf 2}\\
Row -3&{\bf 9}&{\bf 4}\\ 
Row -2&{\bf 16}&{\bf 11}&{\bf 6}&{\bf 1}\\
Row -1&{\bf 23}&{\bf 18}&{\bf 13}&{\bf 8}&{\bf 3}\\
\hline
%\end{tabular}\\
%\begin{tabular}{c c c c c }\hline
%\end{tabular}\\
%\begin{tabular}{c c c c c }
Row 0&{\bf 30}&{\bf 25}&{\bf 20}&{\bf 15}&{\bf 10}&{\bf 5}&{\bf 0}\\
Row 1& 37&32&27&{\bf 22}&{\bf 17}&{\bf 12}&{\bf 7}\\
Row 2&44&39&34&29&24&19&{\bf 14}\\
Row 3& 51&46&41&36&31&26&{\bf 21}\\
Row 4&58&53&48&43&38&33&{\bf 28}\\

....
\end{tabular}
\end{center}
}
\bigskip

We now turn to the case of {\it bar partitions} and the proof of Theorem 4.
Let  $s,t$ be relatively prime odd positive integers. As in \cite{BO}
the $(s,t)$-diagram is divided into 3 parts.
Let $u=\frac{s-1}{2}, v=\frac{t-1}{2}.$ There is a rectangular subdiagram with $u$
rows and $v$ columns with the number $st-s-t$ in its lower left hand
corner and the number $(s+t)/2$ in its upper right hand
corner. We refer to this as the {\it mixed part}. Outside of this 
there are two disjoint subdiagrams. We  refer to the upper
one as the {\it Yin part} and to the lower one as the {\it Yang part}.
In the example below the Yin part is with numbers in {\bf bold} and the
Yang part with numbers in {\it italics}. That two runners are {\it
  conjugate} wrt. $t$ means that the sum of any number on one runner
and any number on the second runner is divisible by $t.$ Similarly we
define conjugate runners wrt. $s.$

The divided (5,7)-{\it diagram}:

{\scriptsize
\begin{center}

\begin{tabular}{c c c c c}
 {\bf 2}\\ 
{\bf 9}&{\bf 4}\\
16&11&6&{\it 1}\\
23&18&13&{\it 8}&{\it 3}\\
%\hline
\end{tabular}
\end{center}
}

In addition we divide the $(s,t)$-abacus into 4 parts: 
Part $A$ consists of the rows with numbers $k < -u.$
Part $B$ is the rows  $-u$ to $-1.$
Part $C$ is the rows  $1$ to $u.$
Part $D$ is the rows $u+1$ to $2u$.
It should be noted, that Part $A$ is the Yin part of the
$(s,t)$-diagram and that Part $B$ contains the Yang part of the
$(s,t)$-diagram.

The rows numbered by the pairs of integers in the following lists 
are called {\it paired}:

$(u+1-j,u+j), j=1,...,u$ are $(C,D)$-paired (placed symmetrically around a
line between rows $u$ and $u+1$)

$(-j,j),j=1,...,u$ are $(B,C)$-paired (placed symmetrically around row
0)

$(-u-j,-(u+1)+j), j=1,...,u $ are $(A,B)$-paired (placed symmetrically around a
line between rows $-u$ and $-(u+1)$.)

\medskip

\noindent {\bf Proof of Theorem 4:} We represent the parts of
$\rho$ as beads on the $(s,t)$-abacus. Since  $\rho$ is a $\bar
t$-core all beads are in the 
top positions on their runners and one of each pair of conjugate
runners wrt. $t$ is empty. Suppose that runner $i$ contains $m_i \ge 0$
beads. We notice that the parts represented by any $s$ consecutive
beads on a runner have different
residue classes modulo $s$ and thus they do not influence the $\bar
s$-core of $\rho.$ 

{\it Step 1:} Remove series of $s$ consecutive beads on runners,
starting from below. 
We still have a $\bar t$-core with the same $\bar s$-core as $\rho.$

After this we assume without loss of generality that $m_i \leq s-1.$

We decompose $m_i=m_i(A)+ m_i(B)+
m_i(C)+m_i(D)+e $ according to the number of beads in the parts $A,B,C,D$
respectively, where $e=0,1$ accounts for a possible bead in row 0. 

{\it Step 2a:} Remove all beads in row 0. Then consider those $i$ for
which $m_i(D)>0.$ 
Do the following:
Remove  $m_i(D) $ pairs of beads from the
$i$'th runner, where each pair of beads is on $(C,D)$-paired rows. 
Modify the $m_i$'s accordingly.
We are then in the situation that $m_i(D)=0$ for all $i$ and we still have a
$\bar t$-core with the same $\bar s$-core as $\rho.$ 

{\it Step 2b:} Consider those $i$ for which $m_i(C)>0.$ Do the
following: Remove  ${\rm min}\{m_i(C),
m_i(B)\} $ pair of beads from the
$i$'th runner  where each pair of beads is on $(B,C)$-paired rows.
Modify the $m_i$'s accordingly. 
We are then in the situation that for each $i$ either $m_i(B)=0$ or $m_i(C)=0.$
We still have a bar partition with the same $\bar s$-core as $\rho.$

{\it Step 2c:} Consider those $i$ for which $m_i(B)>0.$ Do the
following: Remove  ${\rm min}\{m_i(A),
m_i(B)\} $ pair of beads from the
$i$'th runner  where each pair of beads is on $(A,B)$-paired rows.
Modify the $m_i$'s accordingly.
We are then in the situation that for each $i$ at most one of
$m_i(A),m_i(B), m_i(C)$ is nonzero and $m_i(D)=0.$
We still have a bar partition with the same $\bar s$-core as $\rho.$

{\it Step 3:} Any beads left on part C of a runner are moved to part A
of the (empty) conjugate runner. The row number in which a bead is
placed is reduced by $s$. Thus the beads
are moved into the Yin part. 
Any beads left on part B of a runner are moved to part B
of the (empty) conjugate runner, if it is to the right. Thus the beads
are moved into the Yang part. Otherwise 
they are not moved. The row number is unchanged.
Any beads left on part A of a runner are left where they are.

After Step 3 we have reached the bead configuration of a $\bar t$-core
with the same $\bar s$-core as $\rho.$ The only beads left are in
the Yin and Yang parts. 

{\it Step 4:} If possible remove for as long as possible pairs of beads 
in $(A,B)$-paired rows, one bead in Yin, the other in Yang, moving
west to east in both
rows.  Then move the remaining beads as far to the east as possible
in their rows. After Step 4 we have reached the bead configuration of
the $\bar s$-core of  $\rho.$

This is also a $\bar t$-core. All beads are on the top positions on
their runners, due to the fact that all beads are moved as far east 
as possible. Moreover $t$-conjugate runners cannot possibly both
contain beads.  Indeed, such a bead configuration could only occur, when at least one
of the beads has been moved. But the $t$-conjugate runners are placed
symmetrically around a vertical line in the middle of the
$(s,t)$-diagram. Since we start Step 4 in a situation where at least
one of each pair of $t$-conjugate runners is empty, 
this is not possible.  \qed

\medskip

Here is an example illustrating the steps of the proof. 

$(C,D)$-paired rows: (2,3) and (1,4)

$(B,C)$-paired rows: (-1,1) and (-2,2)

$(A,B)$-paired rows: (-3,-2) and (-4,-1) 

{\scriptsize
\begin{center}
\begin{tabular}{l l c c c c c c c }

Runner:& & 2&4&6&1&3&5&0 \\
\hline
\hline
Part A: &Row -4&2\\
& Row -3&9&{\bf 4}\\ 
\hline 
Part B &Row -2&16&{\bf 11}&6&{\bf 1}\\
& Row -1&23&{\bf 18}&13&{\bf 8}&3\\
\hline
%\end{tabular}\\
%\begin{tabular}{c c c c c }
&Row 0&30& 25&20&{\bf 15}&10&{\bf 5}&0\\
\hline
Part C &Row 1&37&32&27&{\bf 22}&17&{\bf 12}&7\\
&Row 2&44&39&34&29&24&{\bf 19}&14\\
\hline
Part D&Row 3&51&46&41&36&31&{\bf 26}&21\\
&Row 4&58&53&48&43&38&33&28\\
%&Row 5&65&60&55&50&45&40&35\\
....
\end{tabular}
\end{center}
}

In step 2a we remove {\bf 5} and {\bf 15} (both in row 0) and the pair
{\bf 19,~26} on
runner 5, rows 2 and 3.
In step 2b we remove the pair the pair {\bf 8,~22} on runner 1, rows -1 and 1. 
In step 2c we remove the pair {\bf 4,~11} on runner 4, rows -3,-2. 
We are left with:

\medskip

{\scriptsize
\begin{center}
\begin{tabular}{l l c c c c c c c }

Runner:& & 2&4&6&1&3&5&0 \\
\hline
\hline
Part A: &Row -4&2\\
& Row -3&9& 4\\ 
\hline 
Part B &Row -2&16& 11&6&{\bf 1}\\
& Row -1&23&{\bf 18}&13& 8&3\\
\hline
%\end{tabular}\\
%\begin{tabular}{c c c c c }
&Row 0&30& 25&20& 15&10& 5&0\\
\hline
Part C &Row 1&37&32&27&22&17&{\bf 12}&7\\
&Row 2&44&39&34&29&24& 19&14\\
%\hline

....
\end{tabular}
\end{center}
}

In step 3 {\bf 12} in the position row 1, runner 5 (Part C) is moved
to {\bf 2} in the position row -4 (Part A). Also {\bf 18} in row -1    
runner 3  is moved to {\bf 3} in the same row on the conjugate runner.

{\scriptsize
\begin{center}
\begin{tabular}{l l c c c c c c c }

Runner:& & 2&4&6&1&3&5&0 \\
\hline
\hline
Part A: &Row -4&{\bf 2}\\
& Row -3&9& 4\\ 
\hline 
Part B &Row -2&16& 11&6&{\bf 1}\\
& Row -1&23& 18&13& 8&{\bf 3} \\
\hline
%\end{tabular}\\
%\begin{tabular}{c c c c c }
&Row 0&30& 25&20& 15&10& 5&0\\
\hline
Part C &Row 1&37&32&27&22&17& 12&7\\
....
\end{tabular}
\end{center}
}

In the final step 4 the pair {\bf 2, 3} is removed, leaving just the
partition (1), a $\bar 7$-core.

\vspace{3ex}

\end{document}